\documentclass[12pt]{amsart}
\usepackage{amscd}
\usepackage{epic,eepic}

\newtheorem{lemma}{Lemma}[section]
\newtheorem{corollary}[lemma]{Corollary}
\newtheorem{theorem}[lemma]{Theorem}
\newtheorem{proposition}[lemma]{Proposition}

\theoremstyle{definition}
\newtheorem{definition}[lemma]{Definition}
\newtheorem{remark}[lemma]{Remark}
\newtheorem{example}[lemma]{Example}

\newtheorem{remark/definition}[lemma]{Remark/Definition}

\def\NZQ{\mathbb}               
\def\NN{{\NZQ N }}

\def\ZZ{{\NZQ Z}}

\def\CC{{\NZQ C}}
\def\opn#1#2{\def#1{\operatorname{#2}}} 

\let\:=\colon
\let\epsilon=\varepsilon
\let\phi=\varphi
\let\kappa=\varkappa

\opn\ini{in}
\def\mm{{\mathfrak m}}
\def\nn{{\mathfrak n}}
\def\pp{{\mathfrak p}}
\def\qq{{\mathfrak q}}
\opn\height{height}

\opn\Ker{Ker} \opn\CoKer{CoKer} \opn\Tor{Tor} \opn\Ext{Ext}
\def\mod{\mathbin{{\textup{mod}}}}
\def\depth{{\operatorname{depth}}}
\opn\projdim{projdim} \opn\reg{reg} \opn\gin{gin} \opn\Cl{Cl}

\let\oldlabel=\label
\def\prellabel{\marginparsep=1em\marginparwidth=80pt
     \def\label##1{\oldlabel{##1}\ifmmode\else\ifinner\else
          \marginpar{{\footnotesize\ \\ \tt
                     ##1}}\fi\fi}}

\begin{document}

\title{Gr\"obner bases, initial ideals and initial algebras}
\author{Winfried Bruns and Aldo Conca}
\address{Universit\"at Osnabr\"uck, FB
Mathematik/Informatik, 49069 Osna\-br\"uck, Germany}
\email{Winfried.Bruns@math.uos.de}
\address{Dipartimento di Matematica, Universit\'a di Genova,
Via Dodeca\-neso 35, 16146 Genova, Italy}
\email{conca@dima.unige.it}

\begin{abstract}
We give an introduction to the theory of initial ideals and
initial algebras with emphasis on the transfer of structural
properties.
\end{abstract}


\maketitle

The notion of Gr\"obner basis of an ideal is the foundation of all
efficient computations in algebraic geometry and commutative
algebra. Highly sophisticated algorithms have been implemented in
several, widely used computer programs.

However, Gr\"obner bases and their analogue for subalgebras are
also important from a purely structural point of view. They allow
us to find deformations of interesting, but ``complicated'' rings
$R$ to simpler objects $R'$ that are defined by monomials and
therefore accessible to combinatorial methods. See \cite{BC} for a
paradigmatic case. In order to transfer the properties that have
been found for $R'$ back to $R$, one has to understand how $R$ and
$R'$ are related. In this article we want to explain this
relationship and to prove some of the basic results about the
passage from $R$ to $R'$, or rather the other way round.

In \cite{BC} we have treated the subject in a similar manner.
However, we hope that some readers will welcome a separate
discussion that is independent from determinantal ideals and
rings. Moreover, the material covered has been slightly expanded
and some proofs are given in more detail.

We are grateful to Tim R\"omer for his careful reading of the
paper and his valuable suggestions.

\section{Initial vector spaces, ideals and subalgebras}
\label{IniVIS}

Let us first recall the definitions and some important properties
of Gr\"obner bases, monomial orders, initial ideals and initial
algebras. For further information on the theory of Gr\"obner bases
we refer the reader to the books by Eisenbud \cite{Eis}, Eisenbud
et al.\ \cite{EGSS}, Greuel and Pfister \cite{GP}, Kreuzer and
Robbiano \cite{KRo}, Sturmfels \cite{Stu2} and Vasconcelos
\cite{Va}. For the so-called Sagbi bases and initial algebras one
should consult Conca, Herzog and Valla \cite{CHV}, Robbiano and
Sweedler \cite{RS}, and \cite[Chapter 11]{Stu2}. Many applications
of Gr\"obner bases are discussed in Buchberger and Winkler
\cite{BW}.

Throughout this section let $K$ be a field, and let $R$ be the
polynomial ring $K[X_1,\dots,X_n]$. A \emph{monomial} (or power
product) of $R$ is an element of the form $X^\alpha=\prod_{i=1}^n
X_i^{\alpha_i}$ with $\alpha\in\NN^n$. A \emph{term} is an element
of the form $\lambda m$ where $\lambda$ is a non-zero element of
$K$ and $m$ is a monomial. Let $M(R)$ be the $K$-basis of $R$
consisting of all the monomials of $R$. Every polynomial $f\in R$
can be written as a sum of terms. The only lack of uniqueness in
this representation is the order of the terms. If we impose a
total order on the set $M(R)$, then the representation is uniquely
determined, once we require that the monomials are written
according to the order, from the largest to the smallest. The set
$M(R)$ is a semigroup (naturally isomorphic to $\NN^n$) and a
total order on the set $M(R)$ is not very useful unless it
respects the semigroup structure.

\begin{definition}\label{monord}
A \emph{monomial order} $\tau$ is a total order
$<_\tau$ on the set $M(R)$ which satisfies the following
conditions:
\begin{itemize}
\item[(a)] $1<_\tau m$ for all the monomials $m\in M(R)\setminus\{1\}$.

\item[(b)] If $m_1,m_2,m_3\in M(R)$ and $m_1<_\tau m_2$, then
$m_1m_3<_\tau m_2m_3$.
\end{itemize}
\end{definition}

From the theoretical as well as from the computational point of
view it is important that descending chains in $M(R)$ terminate:

\begin{remark}\label{wellorder}
A \emph{monomial order} on the set $M(R)$ is a well-order, i.e.
every non-empty subset of $M(R)$ has a minimal element.
Equivalently, there are no infinite descending chains in $M(R)$.

This follows from the fact that every (monomial) ideal in $R$ is
finitely generated. Therefore a subset $N$ of $M(R)$ has only
finitely many elements that are minimal with respect to
divisibility. One of them is the minimal element of $N$.
\end{remark}

We list the most important monomial orders.

\begin{example}\label{lex-revlex}
For monomials $m_1=X_1^{\alpha_1}\cdots X_n^{\alpha_n}$ and
$m_2=X_1^{\beta_1}\cdots\allowbreak X_n^{\beta_n}$ one defines
\begin{itemize}
\item[(a)] the \emph{lexicographic order} (Lex) by
$m_1 <_{\textup{Lex}} m_2$ iff  for some $k$ one has $\alpha_k<
\beta_k$ and $\alpha_i = \beta_i$ for $i<k$;

\item[(b)]  the \emph{degree lexicographic order} (DegLex) by $m_1
<_{\textup{DegLex}} m_2$ iff $\deg(m_1)<\deg(m_2)$ or
$\deg(m_1)=\deg(m_2)$ and $m_1 <_{Lex} m_2$;

\item[(c)] the \emph{(degree) reverse lexicographic order}
(RevLex) by $m_1 <_{\textup{RevLex}} m_2$ iff
$\deg(m_1)<\deg(m_2)$ or $\deg(m_1)=\deg(m_2)$ and for some $k$
one has $\alpha_k> \beta_k$  and $\alpha_i = \beta_i$ for $i>k$.
\end{itemize}
\end{example}

These three monomial orders satisfy $X_1>X_2>\dots>X_n$. More
generally, for every total order on the indeterminates one can
consider the Lex, DegLex and RevLex orders extending the order of
the indeterminates; just change the above definition
correspondingly.

From now on we fix a monomial order $\tau$ on (the monomials of)
$R$. Whenever there is no danger of confusion we will write $<$
instead of $<_\tau$. Every polynomial $f\neq 0$ has an unique
representation
$$
f=\lambda_1 m_1+\lambda_2 m_2 +\dots+\lambda_k m_k
$$
where $\lambda_i \in K\setminus\{0\}$ and $m_1,\dots,m_k$ are
distinct monomials such that $m_1>\dots >m_k$. The \emph{initial
monomial} of $f$ with respect to $\tau$ is denoted by
$\ini_\tau(f)$ and is, by definition, $m_1$.  Clearly one has
$$
\ini_\tau(fg)=\ini_\tau(f)\ini_\tau (g) \eqno{(1)}
$$
and $\ini_\tau(f+g)\leq \max_\tau\{\ini_\tau(f),\ini_\tau(g)\}$.
For example, the initial monomial of the polynomial
$f=X_1+X_2X_4+X_3^2$ with respect to the Lex order is $X_1$, with
respect to DegLex it is $X_2X_4$, and with respect to RevLex it is
$X_3^2$.

Given a $K$-subspace $V\neq 0$ of $R$, we define
$$
M_\tau(V)=\{ \ini_\tau(f) : f\in V\}
$$
and set
$$
\ini_\tau(V)=\mbox{ the $K$-subspace of $R$  generated by }
M_\tau(V).
$$
The space $\ini_\tau(V)$ is called the {\em space of the initial
terms of $V$.}  Whenever there is no danger of confusion we
suppress the reference to the monomial order and use the notation
$\ini(f)$, $M(V)$ and $\ini(V)$.

Any positive integral vector $a=(a_1,\dots, a_n)\in\NN^n$ induces
a graded structure on $R$, called the \emph{$a$-grading}. With
respect to the $a$-grading the indeterminate $X_i$ has degree
$a(X_i)=a_i$. Every monomial $X^\alpha$ is $a$-homogeneous of
$a$-degree $\sum \alpha_ia_i$, and the $a$-degree $a(f)$ of a
non-zero polynomial $f\in R$ is the largest $a$-degree of a
monomial in $f$. Then $R=\bigoplus_{i=0}^\infty R_i$ where $R_i$
is the \emph{$a$-graded component} of $R$ of degree $i$, i.e.\ the
span of the monomials of $a$-degree $i$. With respect to this
decomposition $R$ has the structure of a positively graded
$K$-algebra \cite[Section 1.5]{BH}. The elements of $R_i$ are
\emph{$a$-homogeneous} of $a$-degree $i$. We say that a vector
subspace $V$ of $R$ is \emph{$a$-graded} if it is generated, as a
vector space, by homogeneous elements. This amounts to the
decomposition $V=\bigoplus_{i=0}^\infty V_i$ where $V_i=V\cap
R_i$.

\begin{proposition}\label{basis}
Let $V$ be a $K$-subspace of $R$.
\begin{itemize}
\item[(a)] If $m\in M(V)$ then there exists $f_m\in V$ such that
$\ini(f_m)=m$. The polynomial $f_m$ is uniquely determined if we
further require that the support of $f_m$ intersects $M(V)$
exactly in $m$ and that $f_m$ has leading coefficient $1$.

\item[(b)] $M(V)$ is a $K$-basis of $\ini(V)$.

\item[(c)] The set $\{ f_m : m\in M(V)\}$ is a $K$-basis of $V$.

\item[(d)] If $V$ has finite dimension, then $\dim(V)=\dim(\ini(V))$.

\item[(e)] Let $a\in \NN^n$ be a positive weight vector. Suppose $V$ is
$a$-graded, say $V=\bigoplus_{i=0}^\infty V_i$. Then
$\ini(V)=\bigoplus_{i=0}^\infty \ini(V_i)$. In particular, $V$ and
$\ini(V)$ have the same Hilbert function, i.e.
$\dim(V_i)=\dim(\ini(V)_i)$ for all $i\in\NN$.

\item[(f)] Let $V_1\subseteq V_2$ be $K$-subspaces of $R$. Then
$\ini(V_1)\subseteq \ini(V_2)$ and the (residue classes of the)
elements in $M(V_2)\setminus M(V_1)$ form a $K$-basis of the
quotient space $\ini(V_2)/\ini(V_1)$. Furthermore the set of the
(residue classes of the) $f_m$ with $f_m\in V_2$ and $m \in
M(V_2)\setminus M(V_1)$ is a $K$-basis of $V_2/V_1$ (regardless of
the choice of the $f_m$).

\item[(g)] The set of the (residue classes of the) elements in
$M(R)\setminus M(V)$ is a $K$-basis of $R/V$.

\item[(h)] Let $V_1\subseteq V_2$ be $K$-subspaces of $R$. If
$\ini(V_1)=\ini(V_2)$, then $V_1=V_2$.

\item[(i)] Let $V$ be a $K$-subspace of $R$ and $\sigma,\tau$
monomial orders. If $\ini_\tau(V)\subseteq \ini_\sigma(V)$, then
$\ini_\tau(V)= \ini_\sigma(V)$.
\end{itemize}

\end{proposition}

\begin{proof} (a) and (b) follow easily from the fact that the
monomials form a $K$-basis of $R$. For (a) we have to use that
descending chains in $M(R)$ terminate.

To prove (c) one notes that the $f_m$ are linearly independent
since they have distinct initial monomials. To show that they
generate $V$, we pick any non-zero $f\in V$ and set $m=\ini(f)$.
Then $m\in M(V)$ and we may subtract from $f$ a suitable scalar
multiple of $f_m$, say $g=f-\lambda f_m$, so that
$\ini(g)<\ini(f)$, unless $g=0$. Since $g\in V$, we may repeat the
procedure with $g$ and go on in the same manner. By Remark
\ref{wellorder}, after a finite number of steps we reach $0$, and
$f$ is a linear combination of the polynomials $f_m$ collected in
the subtraction procedure.

(d) and (e) follow from (b) and (c) after the observation that the
element $f_m$ can be taken $a$-homogeneous if $V$ is $a$-graded.

The first two assertions in (f) are easy. For the last we note
that $f_m$ can be chosen in $V_1$ if $m\in\ini(V_1)$.

The residue classes of the $f_m$ with $m \in M(V_2)\setminus
M(V_1)$ are linearly independent modulo $V_1$ since otherwise
there would be a non-trivial linear combination $g=\sum \lambda_m
f_m\in V_1$. But then $\ini(g)\in \ini(V_1)$, a contradiction
since $\ini(g)$ is one of the monomials $m$ which by assumption do
not belong to $M(V_1)$.

To show that the $f_m$ with $m\in M(V_2)\setminus M(V_1)$ generate
$V_2/V_1$ take some non-zero element $f\in V_2$ and set
$m=\ini(f)$. Subtracting a suitable scalar multiple of $f_m$ from
$f$ we obtain a polynomial in $V_2$ with smaller initial monomial
than $f$ (or $0$). If $m\in M(V_1)$, then $f_m\in V_1$. Repeating
the procedure we reach $0$ after finitely many steps. So $f$ can
be written as a linear combination of elements of the form $f_m$
with $m \in M(V_2)\setminus M(V_1)$ and elements of $V_1$, which
is exactly what we want.

(g) is a special case of (f) with $V_2=R$ since in this case we
can take $f_m=m$ for all $m\in M(R)\setminus M(V)$.

(h) follows from (f) since $\ini(V_1)=\ini(V_2)$ implies that the
empty set is a basis of $V_2/V_1$.

Finally, (i) follows from (g) because an inclusion between the two
bases $\{ m\in M(R) : m\not\in M_\tau(V)\}$ and $\{ m\in M(R) :
m\not\in M_\sigma(V)\}$ of the space $R/V$ implies that they are
equal. \end{proof}

\begin{remark/definition}\label{alg-id}
(a) If $A$ is a $K$-subalgebra of $R$, then $\ini(A)$ is also a
$K$-subalgebra of $R$. This follows from equation (1) and from
\ref{basis}(a). The $K$-algebra $\ini(A)$ is called the
\emph{initial algebra of $A$} (with respect to $\tau$).

(b) If $A$ is a $K$-subalgebra of $R$ and $J$ is an ideal of $A$,
then $\ini(J)$ is an ideal of the initial algebra $\ini(A)$. This,
too, follows from equation (1) and from \ref{basis}(a).

(c) If $I$ is an ideal of $R$, then $\ini(I)$ is also an ideal of
$R$. This is a special case of (b) since $\ini(R)=R$.
\end{remark/definition}

\begin{definition}\label{Sag}
Let $A$ be $K$-subalgebra of $R$. A subset $F$ of $A$ is said to
be a \emph{Sagbi basis} of $A$ (with respect to $\tau$) if the
initial algebra $\ini(A)$ is equal to the $K$-algebra generated by
the monomials $\ini(f)$ with $f\in F$.
\end{definition}

If the initial algebra $\ini(A)$ is generated, as a $K$-algebra,
by a set of monomials $G$, then for every $m$ in $G$ we can take a
polynomial $f_m$ in $A$ such that $\ini(f_m)=m$. Therefore $A$ has
a finite Sagbi basis iff $\ini(A)$ is finitely generated. However
it may happen that $A$ is finitely generated, but $\ini(A)$ is
not. The following example is given in \cite{RS} (with a somewhat
different reasoning).

\begin{example}\label{nofinsagbi}
Let $K$ be an arbitrary field, $\tau$ a term order on $K[X,Y]$,
and $A=K[X+Y,XY,XY^2]$. The reader may check that $A$ contains all
monomials $XY^k$, $k\ge 1$. Therefore all these monomials belong
to $\ini(A)$, as well as $X$ if $X>Y$. Now we compute the Hilbert
series of $A$: it is generated by elements of degree $1,2,3$ with
a relation in degree $6$. So
$$
H_A(t)=\frac{1-t^6}{(1-t)(1-t^2)(1-t^3)}=\frac{1-t+t^2}{(1-t)^2}.
$$
Then it is easy to check that the Hilbert function takes the value
$1$ in degree $0$ and $n$ in degree $n$, $n>0$. The algebra
generated by the monomials $XY^k$, $k\ge 0$, has the same Hilbert
function. Since it is contained in the initial algebra (in the
case $X>Y$) it is in fact the initial algebra, but certainly not
finitely generated.

If $Y>X$, one uses that $A$ is also generated by $X+Y$, $YX$, and
$YX^2$.
\end{example}

\begin{definition}\label{Gro}
Let $A$ be a $K$-subalgebra of $R$ and $J$ be an ideal of $A$. A
subset $F$ of $J$ is said to be a \emph{Gr\"obner basis} of $J$
with respect to $\tau$ if the initial ideal $\ini(J)$ is equal to
the ideal of $\ini(A)$ generated by the monomials $\ini(f)$ with
$f\in F$.
\end{definition}

If the initial ideal $\ini(J)$ is generated, as an ideal of
$\ini(A)$, by a set of monomials $G$, then for every $m$ in $G$ we
can take a polynomial $f_m$ in $J$ such that $\ini(f_m)=m$.
Therefore $J$ has a finite Gr\"obner basis iff $\ini(J)$ is
finitely generated. In particular, if $\ini(A)$ is a finitely
generated $K$-algebra, then it is Noetherian and so all the ideals
of $A$ have a finite Gr\"obner basis. Evidently, all the ideals of
$R$ have a finite Gr\"obner basis.

There is an algorithm to determine a Gr\"obner basis of an ideal
of $R$ starting from any (finite) system of generators, the famous
\emph{Buchberger algorithm}. Similarly there is an algorithm that
decides whether a given (finite) set of generators for a
subalgebra $A$ is a Sagbi basis. There also exists a procedure
that completes a system of generators to a Sagbi basis of $A$, but
it does not terminate if the initial algebra is not finitely
generated. If a finite Sagbi basis for an algebra $A$ is known, a
generalization of Buchberger's algorithm finds Gr\"obner bases for
ideals of $A$. We refer the interested readers to the literature
quoted at the beginning of this section.

\section{Initial objects with respect to weights}
\label{IniWeight}

In order to present the deformation theory for initial ideals and
algebras we need to further generalize these notions and consider
initial objects with respect to weights. As pointed out above, any
positive integral weight vector $a=(a_1,\dots, a_n) \in \NN^n$
induces a structure of a positively graded algebra on $R$. Let $t$
be a new variable and set
$$
S=R[t].
$$
For $f=\sum \gamma_i m_i\in R$ with $\gamma_i\in K$ and monomials
$m_i$ one defines the \emph{$a$-homo\-ge\-niza\-tion} $\hom_a(f)$ of
$f$ to be the polynomial
$$
\hom_a(f)=\sum \gamma_i m_i t^{a(f)-a(m_i)}.
$$

Let $a'=(a_1,\dots,a_n,1)\in \NN^{n+1}$.  Clearly, for every $f\in
R$ the element $\hom_a(f)\in S$ is $a'$-homogeneous, and
$f=\hom_a(f)$ iff $f$ is $a$-homogeneous. One has
$$
\begin{aligned}
\ini_a(fg)&=\ini_a(f)\ini_a(g)\\
\hom_a(fg)&= \hom_a(f) \hom_a(g)
\end{aligned}
\qquad\text{for all } f,g\in R. \eqno{(2)}
$$
For every $K$-subspace $V$ of $R$ we set
\begin{align*}
\ini_a(V)&= \text{ the $K$-subspace of $R$ generated by
$\ini_a(f)$ with $f\in V$},\\
\hom_a(V)&=\text{ the $K[t]$-submodule of $S$ generated by $\hom_a(f)$}\\
&\text{\hspace*{2cm}with $f\in V$}.
\end{align*}

If $A$ is a $K$-subalgebra of $R$ and $J$ is an ideal of $A$, then
it follows from (2) that $\ini_a(A)$ is a $K$-subalgebra of $R$
and $\ini_a(J)$ is an ideal of $\ini_a(A)$. Furthermore
$\hom_a(A)$ is a $K[t]$-subalgebra of $S$ and $\hom_a(J)$ is an
ideal of $\hom_a(A)$. As for initial objects with respect to
monomial orders, $\ini_a(A)$ and $\hom_a(A)$ need not be finitely
generated $K$-algebras, even when $A$ is finitely generated. But
if $\ini_a(A)$ is finitely generated, we may find generators of
the form $\ini_a(f_1), \dots, \ini_a(f_k)$ with $f_1,\dots, f_k\in
A$. It is easy to see that the $f_i$ generate $A$. This follows
from the next lemma  in which we use the notation $f^\alpha=\prod
f_i^{\alpha_i}$ for a vector $\alpha \in \NN^k$ and the list
$f=f_1,\dots, f_k$.

\begin{lemma}\label{normalformweight}
Let $A$ be $K$-subalgebra of $R$. Assume that $\ini_a(A)$ is
finitely generated by $\ini_a(f_1), \dots, \ini_a(f_k)$ with
$f_1,\dots, f_k\in A$. Then every $F\in A$ has a representation
$$
F=\sum \lambda_i f^{\beta_i}
$$
where $\lambda_i\in K\setminus\{0\}$ and $a(F)\geq a(f^{\beta_i})$
for all $i$.
\end{lemma}

\begin{proof} By decreasing induction on $a(F)$. The case $a(F)=0$ being
trivial, we assume $a(F)>0$. Since $F\in A$ we have $\ini_a(F)\in
\ini_a(A)=K[\ini_a(f_1),\dots,\allowbreak\ini_a(f_k)]$. Since
$\ini_a(F)$ is an $a$-homogeneous element of the $a$-graded
algebra $\ini_a(A)$, we may write
$$
\ini_a(F)=\sum \lambda_i \ini_a(f^{\alpha_i})
$$
where $a(\ini_a(f^{\alpha_i}))=a(\ini_a(F))$ for all $i$. We set
$F_1=F-\sum \lambda_i f^{\alpha_i}$ and conclude by induction
since $a(F_1)<a(F)$ if $F_1\neq 0$. \end{proof}

The following lemma contains a simple but crucial fact:

\begin{lemma}\label{founda}
Let $A$ be a $K$-subalgebra of $R$ and $J$ be an ideal of $A$.
Assume that $\ini_a(A)$ is finitely generated by $\ini_a(f_1),
\dots, \ini_a(f_k)$ with $f_1,\dots, f_k\in A$. Let
$B=K[Y_1,\dots, Y_k]$ and take presentations
$$
\phi_1:B\to A/J \quad \mbox { and } \quad \phi:B\to
\ini_a(A)/\ini_a(J)
$$
defined by the substitutions $\phi_1(Y_i)=f_i \mod (J)$ and
$\phi(Y_i)=\ini_a(f_i) \mod(\ini_a(J))$. Set $b =(a(f_1), \dots,
a(f_k))\in \NN_+^k$. Then
$$
\ini_b(\Ker \phi_1)=\Ker \phi.
$$

\end{lemma}

\begin{proof}
As a vector space, $\ini_b(\Ker \phi_1)$ is generated by the
elements $\ini_b(p)$ with $p\in \Ker\phi_1$. Set $u=b(p)$. Then we
may write $p=\sum \lambda_i Y^{\alpha_i}+\sum \mu_j Y^{\beta_j}$
where $b(Y^{\alpha_i})=u$ and $b(Y^{\beta_j})<u$. The image
$F=\sum \lambda_i f^{\alpha_i}+\sum \mu_j f^{\beta_j}$ belongs to
$J$, and, hence, $\ini_a(F)\in \ini_a(J)$. Since
$b(Y^\gamma)=a(f^\gamma)$, it follows that $\ini_a(F)=\sum
\lambda_i \ini_a(f^{\alpha_i}$). Thus $\ini_b(p)\in \Ker\phi$, and
this proves the inclusion $\subseteq$.

For the other inclusion we lift $\phi_1$ and $\phi$ to
presentations
$$
\rho_1:B\to A \quad \mbox{  and } \quad \rho:B\to \ini_a(A),
$$
mapping $Y_i$ to $f_i$ and to $\ini_a(f_i)$, respectively. Take a
system of $b$-homogene\-ous generators $G_1$ of the ideal $\Ker
\rho$ of $B$ and a system of $a$-homogeneous generators $G_2$ of
the ideal $\ini_a(J)$ of $\ini_a(A)$. Every $g\in G_2$, being
$a$-homogeneous of degree $u=a(g)$, is of the form $g=\ini_a(g')$,
with $g'\in J$. Then $g'=\sum \gamma_i f^{\alpha_i}+\sum \mu_j
f^{\beta_j}$ with $a(f^{\alpha_i})=u$ and $a(f^{\beta_j})<u$.
Therefore $g=\sum \gamma_i \ini_a(f^{\alpha_i})$.

We choose the canonical preimage of the given representation of
$g$, i.e. $h_g=\sum \gamma_i Y^{\alpha_i}$. Then the set $G_1\cup
\{h_g: g \in G_2\}$ generates the ideal $\Ker \phi$. For all
$g\in G_2$ and $g'$ as above, the canonical preimage of the given
representation of $g'$, i.e. $h=\sum \gamma_i Y^{\alpha_i}+\sum
\mu_j Y^{\beta_j}$ is in $\Ker\phi_1$, and one has
$\ini_b(h)=h_g$.

It remains to show that $g\in \ini_b(\Ker\phi_1)$ for $g\in G_1$.
Every $g\in G_1$ is homogeneous, say of degree $u$, and hence
$g=\sum \lambda_i Y^{\alpha_i}$ with $b(Y^{\alpha_i})=u$. It
follows that $\sum \lambda_i \ini_a(f^{\alpha_i})=0$. Therefore
$\sum \lambda_i f^{\alpha_i}=\sum \mu_j f^{\beta_j}$ with
$a(f^{\beta_j})<u$ by Lemma \ref{normalformweight}. That is,
$g'=\sum \lambda_i Y^{\alpha_i}-\sum \mu_j Y^{\beta_j}$ is in
$\Ker\rho_1$. In particular, $g'\in \Ker\phi_1$ and
$\ini_b(g')=g$.
\end{proof}

A weight vector $a$ and a monomial order $\tau$ on $R$ define a
new monomial order $\tau a$ that ``refines'' the weight $a$ by
$\tau$:
$$
m_1>_{\tau a} m_2 \iff \left\{
\begin{aligned} a(m_1) &> a(m_2) \text{ or }\\
a(m_1)&= a(m_2) \text{ and } m_1 >_{\tau} m_2.
\end{aligned}\right.
$$
We  extend $\tau a$ to $S=R[t]$ by setting:
$$
m_1 t^i >_{\tau a' } m_2 t^j  \iff \left\{
\begin{aligned}
a'(m_1 t^i) &> a'(m_2 t^j ) \mbox{ or }\\
a'(m_1 t^i) &= a'(m_2 t^j )  \mbox{ and } i<j \mbox{ or } \\
a'(m_1 t^i) &= a'(m_2 t^j ) \mbox{ and } i=j \mbox{ and }
m_1>_{\tau} m_2.
\end{aligned}\right.
$$
By construction one has
$$
\ini_{\tau a }(f)=\ini_{\tau a' }(\hom_a(f))\quad \mbox{for all }
f\in R,\ f\neq 0.
$$
Given a $K$-subspace $V$ of $R$, we let $VK[t]$ denote the
$K[t]$-submodule of $S$ generated by the elements in $V$.

\begin{proposition}\label{deformation}
Let $a\in\NN^n$ be a positive integral vector and $\tau$ be a
monomial order on $R$. For every $K$-subspace $V$  of $R$ one has:
\begin{itemize}
\item[(a)] $\ini_{ \tau a}(V)=\ini_{ \tau a}(\ini_
a(V))=\ini_{\tau}(\ini_a(V))$,
\item[(b)] If either $\ini_\tau(V)\subseteq \ini_ a(V)$ or
$\ini_\tau(V)\supseteq \ini_ a(V)$, then $\ini_\tau(V)=\ini_
a(V)$,
\item[(c)] $\ini_{ \tau a}(V)K[t]=\ini_{ \tau a'}(\hom_a(V))$,
\item[(d)] The quotient $S/\hom_a(V)$ is a free $K[t]$-module.
\end{itemize}

\end{proposition}

\begin{proof} (a) Note that
$\ini_{\tau a}(f)=\ini_{\tau a}(\ini_a(f))=\ini_{\tau}(\ini_a(f))$
holds for every $f\in R$. It follows that the first space is
contained in the second and in the third. On the other hand, since
$\ini_a(V)$ is $a$-homogeneous, the monomials in its initial space
are initial monomials of $a$-homogeneous elements. But every
$a$-homogeneous element in $\ini_a(V)$ is of the form $\ini_a(f)$
with $f\in V$. This gives the other inclusions.

(b) If one of the two inclusions holds, then an application of
$\ini_\tau(..)$ to both sides yields that $\ini_\tau(V)$ either
contains or is contained in $\ini_\tau(\ini_a(V))$. By (a) the
latter is $\ini_{ \tau a}( V)$. Then by Proposition \ref{basis}(i)
we have that $\ini_\tau(V)=\ini_{ \tau a}( V)$. Next we may apply
\ref{basis}(h) and conclude that $\ini_\tau(V)=\ini_ a(V)$.

(c) For every $f\in R$ one has $\ini_{\tau a'}(\hom_a(f))
=\ini_{\tau a}(f)$. Thus $\ini_{ \tau a}(V)K[t]\allowbreak
\subseteq \ini_{ \tau a'}(\hom_a(V))$. On the other hand,
$\hom_a(V)$ is an $a'$-homogeneous space. Therefore its initial
space is generated by the initial monomials of its
$a'$-homogeneous elements. An $a'$-homogeneous element of degree,
say, $u$ in $\hom_a(V)$ has the form $g=\sum_{i=1}^k \lambda_i
t^{\alpha_i} \hom_a(f_i)$ where $f_i\in V$ and
$\alpha_i+a(f_i)=u$. If $\alpha_i=\alpha_j$ then $a(f_i)=a(f_j)$
and $\hom_a(f_i+f_j)=\hom_a(f_i)+\hom_a(f_j)$. In other words, we
may assume that the $\alpha_i$ are all distinct and, after
reordering if necessary, that $\alpha_i<\alpha_{i+1}$. Then
$\ini_{\tau a'}(g)= t^{\alpha_1} \ini_{\tau a'}(\hom(f_1))=
t^{\alpha_1} \ini_{\tau a}(f_1)$. This proves the other inclusion.

(d) By (c) and Proposition \ref{basis}(b) the (classes of the)
elements $t^\alpha m$, $\alpha\in \NN$, $m\in M(R)\setminus M(V)$,
form a $K$-basis of $S/\hom_a(V)$. This implies that the set
$M(R)\setminus M(V)$ is a $K[t]$-basis of $S/\hom_a(V)$.
\end{proof}

The next proposition connects the structure of $R/I$ with that of
$R/\ini_a(R)$:

\begin{proposition}\label{flatfamily}
For every ideal $I$ of $R$ the ring $S/\hom_a(I)$ is a free
$K[t]$-module. In particular $t-\alpha$ is a non-zero divisor on
$S/\hom_a(I)$ for every $\alpha\in K$. Furthermore
$S/(\hom_a(I)+(t))\cong R/\ini_a(I)$ and
$S/(\hom_a(I)+(t-\alpha))\cong R/I$ for all $\alpha\neq 0$.
\end{proposition}

\begin{proof}
The first assertion follows from \ref{deformation}(d). It implies
that every non-zero element of $K[t]$ is a non-zero divisor on
$S/\hom_a(I)$. For $S/(\hom_a(I)\allowbreak+(t))\cong R/\ini_a(I)$
it is enough that $\hom_a(I)+(t)=\ini_a(I)+(t)$. This is easily
seen since for every $f\in R$ the polynomials $\ini_a(f)$ and
$\hom_a(f)$ differ only by a multiple of $t$. To prove that
$S/(\hom_a(I)+(t-\alpha))\cong R/I$ for every $\alpha\neq 0$, we
consider the graded isomorphism $\psi: R\to R$ induced by
$\psi(X_i)=\alpha^{-a_i}X_i$. One checks that
$\psi(m)=\alpha^{-a(m)}m$ for every monomial $m$ of $R$ and that
$\hom_a(f)-\alpha^{a(f)} \psi(f)$ is a multiple of $t-\alpha$ for
all the $f\in R$. So $\hom_a(I)+(t-\alpha)= \psi(I)+(t-\alpha)$,
which implies the desired isomorphism. \end{proof}

\section{The transfer of arithmetic and homological properties}
\label{Transfer}

Now we use Proposition \ref{flatfamily} for comparing $R/I$ with
$R/\ini_a(I)$.

\begin{proposition}\label{transfer}
\begin{itemize}
\item[(a)] $R/I$ and $R/\ini_a(I)$ have the same Krull dimension.
\item[(b)] The following properties are passed from $R/\ini_a(I)$
on to $R/I$: being reduced, a domain, a normal domain,
Cohen-Macaulay, Gorenstein.
\item[(c)] Suppose that $I$ is graded with respect to some positive
weight vector $b$. Then $\ini_a(I)$ is $b$-graded, too, and the
Hilbert functions of $R/I$ and $R/\ini_a(I)$ coincide.
\end{itemize}
\end{proposition}

\begin{proof}
Let us start with (b). The bridge between $R/I$ and $R/\ini(I)$ is
formed by $A=S/\hom_a(I)$. We have representations $R/I=A/(t-1)$
and $R/\ini(I)=A/(t)$. So we must show that the properties under
consideration first \emph{ascend} from the residue class ring
$A/(t)$ to $A$ and then \emph{descend} from $A$ to $A/(t-1)$.
\smallskip

\noindent\emph{Ascent.}\enspace The $K$-algebra $A$ is positively
graded. Let $\mm$ denote its maximal ideal generated by the
residue classes of the indeterminates. Set $A'=A_\mm$. Then we
have the following commutative diagram in which all maps are the
natural ones:
$$
\begin{CD}
A@>>> A'\\
@VVV @VVV\\
A/(t)@>>>A'/(t).
\end{CD}
$$

(i) We start from $A/(t)$. The passage to its localization
$A'/(t)$ with respect to the maximal ideal $\mm/(t)$ preserves all
the properties under consideration.

(ii) Now we have to go up from $A'/(t)$ to $A'$ itself. It is
elementary to show that $A'$ is reduced or an integral domain if
$A'/(t)$ has this property (see \cite[Proof of 2.2.3]{BH} for the
prototype of such an argument). Normality is covered by the next
lemma. For the Cohen-Macaulay and Gorenstein property the
conclusion is contained in \cite[2.1.3 and 3.1.9]{BH}.

(iii) Finally, one observes that $A$ has one of the properties
mentioned if and only if its localization $A'=A_\mm$ does so. In
fact, all of the properties depend only on the localizations of
$A$ with respect to graded prime ideals, and such localizations
are localizations of $A'$ (see \cite[Section 1.5 and Chapters 2
and 3]{BH}, in particular \cite[2.1.27, 3.6.20]{BH}). For the only
non-local property, namely that of being an integral domain, one
notes that $\mm$ contains all the associated prime ideals of
$A$.\smallskip

\noindent\emph{Descent.}\enspace It remains to transfer the
properties in (b) to $A''=A/{(t-1)}\allowbreak \cong R/I$. At this point one
should observe that $A''$ is not merely a residue class ring
modulo a non-zero-divisor, but in fact the dehomogenization of $A$
with respect to the degree $1$ element $t$. So $A''$ is the degree
$0$ component of the graded ring $A[t^{-1}]$, and $A[t^{-1}]$ is
just the Laurent polynomial ring in the variable $t$ over $A''$.
(This is not hard to see; cf.\ \cite[Section 1.5]{BH}. The main
point is that the surjection $A\to A''$ factors through
$A[t^{-1}]$ and that the latter ring has a homogeneous unit of
degree $1$.) Finally, each of the properties descends from the
Laurent polynomial ring to $A''$. The proof of (b) is complete.

For (a) one follows the same chain of descents and ascents:
\begin{align*}
\dim R/I &=\dim A''=\dim A''[t,t^{-1}]-1=\dim A[t^{-1}]-1\\
&=\dim A-1=\dim R/\ini(A).
\end{align*}
For the equation $\dim A=\dim A[t^{-1}]$ one has to use that $t$
is a non-zero-divisor in the affine $K$-algebra $A$: it can not be
contained in all maximal ideals $\nn$ of $A$ for which $\dim
A=\dim A_\nn$. In fact, let $\pp$ be a minimal prime ideal of $A$
with $\dim A=\dim A/\pp$. Then all maximal ideals $\nn\supset\pp$
have $\dim A_\nn=\dim A$, and $\pp$ is their intersection. But
$t\notin\pp$. For the very last equation one can use that $t$ is a
homogeneous non-zero-divisor in the positively graded ring $A$.

For (c) one first notes that $\ini_a(I)$ is $b$-graded, since the
initial form of a $b$-homogeneous element is $b$-homogeneous, too.
We refine the weight $a$ by a monomial order $\tau$ and derive the
chain of equations
$$
H(R/\ini_a(I))=H(R/\ini_\tau(\ini_a(I)))=H(R/\ini_{\tau a}(I))
=H(R/I)
$$
for the Hilbert function $H(..)$ from \ref{basis}(e) and
\ref{deformation}(a). \end{proof}

We have to add a lemma already used in the proof above. The local
rings used there are catenary as is every localization of an
affine $K$-algebra.

\begin{lemma}\label{normup}
Let $A$ be a catenary noetherian local ring and $t$ a
non-zero-divisor of $A$. If $A/(t)$ is normal, then so is $A$.
\end{lemma}

\begin{proof}
We must show that $A$ has the Serre properties $(R_1)$ and $(S_2)$
if these hold for $A/(t)$. Let $\pp$ be a prime ideal of $A$ with
$\height\pp\le 1$. If $t\in\pp$, then $\overline\pp=\pp/(t)$ is a
minimal prime ideal of $A/(t)$, and the regularity of
$(A/(t))_{\overline\pp}=A_\pp/(t)$ implies that of $A_\pp$. If
$t\notin\pp$, we choose a minimal prime overideal $\qq$ of
$\pp+(t)$. Since $A$ is catenary, we must have $\height
\qq=\height\pp+1$. Moreover, $\height
\qq/(t)=\height\qq-1=\height\pp$. It follows that
$(A/(t))_{\overline\qq}$ is regular. So $A_\qq$ and its
localization $A_\pp$ are regular.

Suppose now that $\height\pp\ge 2$. We must show that $\depth
A_\pp\ge 2$. If $t\in\pp$, then we certainly have $\depth
(A/(t))_{\overline\pp}\ge 1$, since $(A/(t))_{\overline\pp}$ is
regular or has depth at least $2$. Otherwise we take $\qq$ as
above. Then $\depth (A/(t))_{\overline\qq}\ge 2$, and $\depth
A_\qq\ge 3$. We choose $u\neq 0$ in $\pp$. If $\depth A_\pp=1$,
then $\pp/(u)$ is an associated prime ideal of $A/(u)$. Moreover,
we have $\depth A_\qq/(u)\ge 2$, and $\dim A_\qq/\pp A_\qq=1$.
This is a contradiction to \cite[1.2.13]{BH}: for a local ring $R$
one has $\depth R\le \dim R/\pp$ for all associated prime ideals
$\pp$ of $R$.
\end{proof}

Very often one wants to compare finer invariants of $R/\ini_a(I)$
and $R/I$, for example if $I$ is a graded ideal of $R$ with
respect to some other weight vector $b$. The next proposition
shows that the comparison is possible for graded components of
Tor-modules. The vector space dimensions in the proposition are
called \emph{graded Betti numbers}.

\begin{proposition}\label{flatfam}
Let $a,b$ positive integral vectors and let $J,J_1,J_2$ be
$b$-ho\-mo\-ge\-neous ideals of $R$ with $J \subseteq J_1$ and
$J\subseteq J_2$. Then $\ini_a(J), \ini_a(J_1), \ini_a(J_2)$ are
also $b$-homogeneous ideals, and one has
$$
\dim_K \Tor^{R/J}_i(R/J_1, R/J_2)_j\leq \dim_K
\Tor^{R/\ini_a(J)}_i ( R/\ini_a(J_1), R/\ini_a(J_2))_j
$$
where the graded structure on the $\Tor$-modules is inherited from
the $b$-graded structure of their arguments.
\end{proposition}

\begin{proof}
On $S$ we introduce a bigraded structure, setting $\deg
X_i=(b_i,a_i)$ and $\deg t=(0,1)$. The ideals $I=\hom_a(J)$,
$I_1=\hom_a(J_1)$ and $I_2=\hom_a(J_2)$ are then bigraded and so
are the algebras they define. We need a standard result in
homological algebra: if $A$ is a ring, $M,N$ are $A$-modules and
$x$ is a non-zero-divisor on $A$ as well as on $M$ then
$\Tor^A_i(M,N/xN)\cong \Tor^{A/xA}_i(M/xM, N/xN)$. (It is
difficult to find an explicit reference; for example, one can use
\cite[1.1.5]{BH}.) If, in addition, $x$ is a non-zero-divisor also
on $N$, then we have the short exact sequence $0\to N\to N\to
N/xN\to 0$. It yields the exact sequence
$$
0\to \CoKer \phi_i  \to \Tor^{A/xA}_i(M/xM, N/xN)\to \Ker
\phi_{i-1}\to 0
$$
where $\phi_i$ is multiplication by $x$ on $\Tor^A_i(M,N)$.

Set $A=S/\hom_a(J)$, $M=S/\hom_a(J_1)$, $N=S/\hom_a(J_2)$ and
$T_i=\Tor^{A}_i(M,N)$. Since the modules involved are bigraded,
so is $T_i$. Let $T_{ij}$ be the direct sum of all the components
of $T_i$ of bidegree $(j,k)$ as $k$ varies. Since $T_i$ is a
finitely generated bigraded $S$-module, $T_{ij}$ is a finitely
generated and graded $K[t]$-module (with respect to the standard
grading of $K[t]$). So we may decompose it as
$$
T_{ij}=F_{ij} \oplus G_{ij}
$$
where $F_{ij}$ is the free part and $G_{ij}$ is the torsion part,
which, being $K[t]$-graded, is a direct sum of modules of the form
$K[t]/(t^a)$ for various $a>0$. Denote the minimal number of
generators of $F_{ij}$ and $G_{ij}$ as $K[t]$-modules by $f_{ij}$
and $g_{ij}$, respectively. Now we consider the $b$-homogeneous
component of degree $j$ of the above short exact sequence with
$x=t$, which is a non-zero-divisor by Proposition
\ref{deformation}(d). It follows that
$$
\dim_K \Tor^{R/\ini_a(J)}_i ( R/\ini_a(J_1),
R/\ini_a(J_2))_j=f_{ij}+g_{ij}+g_{i-1,j}.
$$
If we take $x=t-1$ instead of $x$, then we have
$$
\dim_K \Tor^{R/J}_i ( R/J_1, R/J_2)_j=f_{ij}
$$
and this shows the desired inequality. \end{proof}

\begin{remark}
One can prove an analogous inequality for Ext-modules. However,
some care is advisable: the homological degree $i$ changes to
$i-1$ when one passes from $A$ to the residue class rings modulo
$t$ and $t-1$ (Lemma of Rees \cite[3.1.16]{BH}).
\end{remark}

Note that one can use Proposition \ref{flatfam} to transfer the
Cohen-Macaul\-ay and Gorenstein properties from $R/\ini_a(I)$ to
$R/I$ if $I$ is $b$-graded. We content ourselves with a comparison
of two important invariants:

\begin{corollary}\label{pdreg}
Under the hypotheses of \ref{flatfam} one has
$$
\projdim_R R/I\le \projdim_R R/\ini_a(I).
$$
If $(a=(1,\dots,1)$, then
$$
\reg_R R/I\le \reg_R R/\ini(I).
$$
\end{corollary}

\begin{proof}
For both invariants this is an immediate consequence of the
proposition, for the projective dimension
$$
\projdim_R R/I=\max\{i:\Tor_i^R(R/I,K)\neq 0\}
$$
as well as for the Castelnuovo-Mumford regularity
$$
\reg_R R/I=\max\{j-i:\Tor_i^R(R/I,K)_j\neq 0\}.
$$
(In its definition one assumes that all indeterminates have degree
$1$.)
\end{proof}

\begin{remark}\label{gin}
As we will see in Proposition \ref{initbyweight} every monomial
order $\tau$ can be approximated by a weight vector $a$, as long
as one only wants to compute the initial ideals of finitely many
ideals. Therefore Corollary \ref{pdreg} applies also to initial
ideals defined by monomial orders.

While the inequalities in the previous corollary are strict in
general, they turn into equalities in an important special case,
namely when $\tau$ is the RevLex order, and the initial ideal is
formed after a generic linear transformation $\gamma$ of the
coordinates. Then $\gin(I)=\ini_\tau(\gamma(I))$ is called the
\emph{generic initial ideal}. One has $\projdim_R R/I=\projdim_R
R/\gin(I)$ and $\reg_R R/I=\reg_R R/\gin(I)$; see \cite[19.11 and
20.21]{Eis} for this theorem of Bayer and Stillman. For further
results comparing single Betti numbers of $R/I$ and $R/\gin(I)$
see Bayer, Charalambous and Popescu \cite{BCP} and Aramova, Herzog
and Hibi \cite{AHH}.
\end{remark}

If $I$ is graded with respect to the ordinary weight $(1,\dots,1)$
then it makes sense to ask for the Koszul property of $R/I$. By
definition, $R/I$ is Koszul if $\Tor_i^{R/I}(R/\mm,\allowbreak
R/\mm)_j$ is non-zero only for $i=j$. Backelin and Fr\"oberg
\cite{BF} give a detailed discussion of this class of rings.

\begin{corollary}\label{Koszul}
Suppose that $I$ is a graded ideal with respect to the weight
$(1,\dots,1)$.
\begin{itemize}
\item[(a)] If $R/\ini_a(I)$ is Koszul for some positive weight
$a$, then $R/I$ is Koszul.
\item[(b)] In particular, if $\ini_a(I)$ is generated by
degree $2$ monomials, then $R/I$ is Koszul.
\end{itemize}
\end{corollary}

\begin{proof}
(a) follows directly from Proposition \ref{flatfam}. For (b) one
uses a theorem of Fr\"oberg \cite{Fr}: if $J$ is an ideal
generated by quadratic monomials, then the algebra $R/J$ is
Koszul.
\end{proof}

In order to apply the previous results to initial objects defined
by monomial orders we have to approximate such orders by weight
vectors. This is indeed possible, provided only finitely many
monomials have to be considered.

\begin{proposition}\label{initbyweight}
Let $\tau$ be a monomial order on $R$.
\begin{itemize}
\item[(a)] Let $\{(m_1, n_1),\dots, (m_k,n_k)\}$ be a finite
set of pairs of monomials such that $m_i>_\tau n_i$ for all $i$.
Then there exists a positive integral weight $a\in \NN_+^n$ such
that $a(m_i)>a(n_i)$ for all $i$.

\item[(b)] Let $A$ be a $K$-subalgebra  of $R$ and $I_1,\dots,I_h$ be ideals
of $A$. Assume that $\ini_\tau(A)$ is finitely generated as a
$K$-algebra. Then there exists a positive integral weight $a\in
\NN_+^n$ such that $\ini_\tau(A)=\ini_a(A)$ and
$\ini_\tau(I_i)=\ini_a(I_i)$ for all $i=1,\dots, h$.
\end{itemize}
\end{proposition}

\begin{proof} (a) Set $m_i=X^{\alpha_i}$ and $n_i=X^{\beta_i}$ and
$\gamma_i=\alpha_i-\beta_i\in \ZZ^n$. Let $\Gamma$ be the $k\times
n$ integral matrix whose rows are the vectors $\gamma_i$. We are
looking for a positive column vector $a$ such that the
coefficients of the vector $\Gamma a$ are all $>0$. Suppose, by
contradiction, there is no such $a$. Then (one version of the
famous)  Farkas Lemma (see Schrijver \cite[Section 7.3]{Sij}) says
that there exists a linear combination $v=\sum c_i \gamma_i$ with
non-negative integral coefficients $c_i\in \NN$  such that $v\leq
0$, that is $v=(v_1,\dots, v_n)$ with $v_i\leq 0$. Then it follows
that $\prod_i m_i^{c_i}X^{-v}=\prod_i n_i^{c_i}$, which
contradicts our assumptions because the monomial order is
compatible with the semigroup structure.

(b) Let $F_0$ be a finite Sagbi basis of $A$, let $F_i$ be a
finite Gr\"obner basis of $I_i$ and set $F= \bigcup_i F_i$.
Consider the set $U$ of pairs of monomials $(\ini(f),m)$ where
$f\in F$ and $m$ is any non-initial monomial of $f$. Since $U$ is
finite, by (a) there exists $a\in \NN_+^n$ such that
$\ini_a(f)=\ini_\tau(f)$ for every $f\in H$. We show $a$ has the
desired property. Set $V_0=A$ and $V_i=I_i$. By construction the
(algebra for $i=0$ and ideal for $i>0$) generators of the
$\ini_\tau(V_i)$ belong to $\ini_a(V_i)$ so that
$\ini_\tau(V_i)\subseteq \ini_a(V_i)$. But then, by Proposition
\ref{deformation}(b), we may conclude that
$\ini_\tau(V_i)=\ini_a(V_i)$. \end{proof}

The main theorem of this section summarizes what we can say about
the transfer of ring-theoretic properties from initial objects.
For the Koszul property of subalgebras we must allow a
``normalization'' of degree. Suppose that $b$ is a positive weight
vector $b$, and suppose that a subalgebra $A$ is generated by
elements $f_1,\dots,f_s$ of the same $b$-degree $e\in\NN$. Then
every element $g$ of $A$ has $b$-degree divisible by $e$, and
dividing the $b$-degree by $e$ we obtain the \emph{$e$-normalized
$b$-degree} of $g$.

\begin{theorem}\label{reprbywei}
Let $\ini(..)$ denote the initial objects with respect to a
positive integral vector $a\in \NN^n$ or to a monomial order
$\tau$ on $R$. Let $A$ be a $K$-subalgebra of $R$ and $J$ be an
ideal of $A$. Suppose that $\ini(A)$ is finitely generated.
\begin{itemize}
\item[(a)] One has $\dim A/J=\dim \ini(A)/\ini(J)$.
\item[(b)]  If $\ini(A)/\ini(J)$ is reduced, a domain, a
normal domain, Cohen-Macaulay, or Gorenstein, then so is $A/J$.
\item[(c)] Let $b$ be a positive weight vector, and suppose that
$A$ and $J$ are $b$-graded. Then $A/J$ and $\ini(A)/\ini(J)$ have
the same Hilbert function.
\item[(d)] If, in addition to the hypothesis of (c), $\ini(A)/\ini(J)$
is Koszul with respect to $e$-normalized $b$-degree for some $e$,
then so is $A/J$.
\end{itemize}
\end{theorem}

\begin{proof} If the initial objects are formed with respect to a
monomial order then, by \ref{initbyweight}, we may represent them
as initial objects with respect to a suitable  positive integral
weight vector. Therefore in both cases the initial objects are
taken with respect to a positive integral weight $a$. By Lemma
\ref{founda} there exist a polynomial ring, say $B$, an ideal $H$,
and a positive weight $c$ such that $B/H\cong A/J$ and
$B/\ini_c(H)\cong \ini(A)/\ini(J)$. Furthermore, under the
hypothesis of (c), the weight $b$ can be lifted from the
generators of $\ini(A)$ to the indeterminates of $B$. Now the
theorem follows from Proposition \ref{transfer} and Lemma
\ref{Koszul}. \end{proof}

The theorem is usually applied in two extreme cases. In the first
case $A=R$, so that $\ini(A)=R$, and in the second case $H=0$, so
that  $\ini(J)=0$.

There is a special instance of the theorem that deserves a
separate statement.

\begin{corollary}\label{mainc2}
Let $A$ be $K$-subalgebra of $R$, and suppose that $\ini(A)$ is
generated by finitely many monomials (e.g.\ if it is finitely
generated and the initial algebra is taken with respect to a
monomial order). If $\ini(A)/\ini(I)$ is normal, then $A/I$ is
normal and Cohen-Macaulay.
\end{corollary}

\begin{proof}
The hypothesis implies that $\ini(I)$ is a prime ideal in the
affine semigroup ring $\ini(A)$. But then the natural homomorphism
$\ini(A)\to\ini(A)/\ini(I)$ splits as a ring homomorphism (see
\cite[Section 6.1]{BH}). It follows that $\ini(A)/\ini(I)$ is
itself a normal affine semigroup ring. By a theorem of Hochster
\cite[6.3.5]{BH} such a ring is Cohen-Macaulay.
\end{proof}

Sometimes one of the implications in Theorem \ref{reprbywei} can
be reversed:

\begin{corollary}\label{GorIni}
Let $b$ be a positive weight vector, and suppose that the
$K$-subalgebra $A$ is $b$-graded and has a Cohen-Macaulay initial
algebra $\ini(A)$. Then $A$ is Gorenstein iff $\ini(A)$ is
Gorenstein.
\end{corollary}

\begin{proof}
Since $\ini(A)$ is Cohen-Macaulay, $A$ is Cohen-Macaulay as well.
So both algebras are positively graded Cohen-Macaulay domains. By
a theorem of Stanley \cite[4.4.6]{BH}, the Gorenstein property of
such rings depends only on their Hilbert function, and both
algebras have the same Hilbert function. \end{proof}

We want to extend Theorem \ref{reprbywei} in such a way that it
allows us to determine the canonical module of $A/I$. First a
lemma that covers the most difficult step in the passage from
$\ini(A)/\ini(I)$ to $A/I$.

\begin{lemma}\label{cano-up}
Let $R$ be a positively graded algebra over a field $K$ and $C$ a
finitely generated graded $R$-module. Suppose that $t\in R$ is a
homogeneous non-zero-divisor for both $R$ and $C$. Then $C$ is the
canonical module of $R$ (up to the shift $a$) if (and only if)
$C/tC$ is the canonical module of $R/(t)$ (up the shift $a+\deg
t$).
\end{lemma}

\begin{proof}
Let $d=\dim R$. Then $\dim R/(t)=d-1$. We can assume that $a=0$,
shifting $C$ and $C/tC$ by $-a$ if necessary. By the lemma of Rees
(for example, see \cite[3.1.16 and 4.2.40]{BH}) we have
$$
\Ext_{R/(t)}^i\bigl(K,(C/tC)(\deg t)\bigr)\cong\Ext_R^{i+1}(K,C)=
\begin{cases}0& i\neq d-1,\\K&i=d-1.\end{cases}
$$
(with $K$ in degree $0$). This property is exactly the definition
of the graded canonical module; see \cite[Section 3.6]{BH}.
\end{proof}

\begin{theorem}\label{CanIni}
Let $A$ be a subalgebra of $R$ as in Theorem \ref{reprbywei}, and
$I\subseteq J$ ideals of $A$. Suppose that $\ini(A)/\ini(I)$ and,
hence, $A/I$ are Cohen-Macaulay.
\begin{itemize}
\item[(a)] If $\ini(J)/\ini(I)$ is the canonical module of
$\ini(A)/\ini(I)$, then $J/I$ is the canonical module of $A/I$.
\item[(b)] Suppose in addition that $A,I,J$ are $b$-graded with
respect to a positive weight and $\ini(J)/\ini(I)$ is the
canonical module of $\ini(A)/\ini(I)$ (up to a shift). Then $J/I$
is the graded canonical module (up to the same shift).
\end{itemize}
\end{theorem}

\begin{proof}
(a) As in the proof of Theorem \ref{reprbywei} we may assume that
the initial objects are defined by a weight vector. Then we choose
representations $A/I\cong B/I_1$, $A/J\cong B/I_2$,
$\ini(A)/\ini(I)\cong B/\ini(I_1)$, $\ini(A)/\ini(J)\cong
B/\ini(I_2)$ as in Lemma \ref{founda}. This reduces the problem to
the situation of Proposition \ref{transfer}: $R$ is a polynomial
ring over $K$, $I\subset J$ are ideals, and $R/\ini(I)$ is
Cohen-Macaulay with canonical module $\ini(J)/\ini(I)$.

Again one passes to the homogenized objects in $S=R[t]$. Note that
$t$ and $t-1$ are non-zero-divisors modulo $\hom(I)$ and
$\hom(J)$. Set $\bar S= S/\hom(I)$ and $\bar J=\hom(J)/\hom(I)$.
Then $\bar S/t\bar S\cong R/\ini(I)$ and
$$
\bar J/t\bar J\cong (\bar J+t\bar S)/(t\bar S)\cong
\ini(J)/\ini(I).
$$
By Lemma \ref{cano-up} we therefore conclude that $\bar J$ is the
canonical module of $\bar S$. But we also have $\bar J/(t-1)\bar
J\cong (\bar J+(t-1)\bar S)/(t-1)\bar S\cong J/I$. This shows that
$J/I$ is the canonical module of $R/I$.

(b) It only remains to control the shift. This can be done via the
Hilbert functions of $A/I$ and $J/I$ on the one side and those of
$\ini(A)/\ini(I)$ and $\ini(J)/\ini(I)$ on the other (see
\cite[4.4.5]{BH}). But the Hilbert functions of the objects
corresponding to each other via $\ini$ coincide, and the claim
follows.
\end{proof}

\begin{remark} In addition to Cohen-Macaulay and Gorenstein
rings one can also consider those with rational singularities (in
characteristic $0$) or $F$-rational singularities (in
characteristic $p$). They behave well under the deformation to the
initial objects. See \cite{CHV} for a more detailed discussion.

In particular $A/I$ is ($F$-)rational under the hypotheses of
Corollary \ref{mainc2}.
\end{remark}

\begin{remark}\label{worse}
All the results above suggest that the numerical invariants and
structural properties can only improve in the direction from
$B'=\ini(A)/\ini(I)$ to $B=A/I$. However, some caution is
advisable.

(a) If both algebras are normal Noetherian domains, then one can
consider their divisor class groups $\Cl(B)$ and $\Cl(B')$. A
potential theorem comparing $\Cl(B)$ and $\Cl(B')$ could be that
$\Cl(B)$ is always of the form $G/H$ where $H\subset
G\subset\Cl(B')$. This is not the case as the following example
indicates.

Choose $A=\CC[X^2-Z^2, XY, Y^2, YZ]\cong
\CC[T.U.V.W]/(U^2-TV-W^2]$, its initial algebra with respect to
Lex is $\ini(A)=\CC[X^2, XY, Y^2, YZ]\cong \CC[T,U,V,W]/(U^2-TV)$.
For the verification of the claims in this statement it is enough
to observe that $\CC[X^2, XY, Y^2, YZ]$ has indeed the
representation given, and that $U^2-TV-W^2$ is a relation of the
generators of $A$. The rest follows from Hilbert function
arguments.

According to Fossum \cite[11.4]{F} $A$ has divisor class group
isomorphic to $\ZZ$ since the quadratic form defining it is
non-degenerate (and $\CC$ is algebraically closed). However,
$\ini(A)$ has class group $\ZZ/2\ZZ$.

(b) Another (and related invariant) is the Grothendieck group
$K_0(R)$. Gubeladze \cite{Gu} has given an example of an algebra
$A$ with $K_0(A)\neq 0$ for which $\ini(A)$ has trivial $K_0$.
\end{remark}

\end{document}